\documentclass[11pt]{article}
\usepackage[dvips]{graphicx}
\usepackage{amssymb}
\usepackage{color}
\usepackage{amsmath}
\usepackage{amsthm}
\usepackage{latexsym}
\usepackage[english]{babel}
\usepackage{hyperref}
\newcommand{\R}{\mathbb{R}}
\newcommand{\Z}{\mathbb{Z}}

\newcommand{\D}{\mathsf{D}}

\newtheorem{theorem}{Theorem}[section]

\newtheorem{lemma}[theorem]{Lemma}

\newtheorem{proposition}[theorem]{Proposition}

\begin{document}
\title{Super-Symmetric Coupling: Existence and Multiplicity}
\author{Ali Maalaoui$^{1}$}
\addtocounter{footnote}{1}
\footnotetext{Department of mathematics and natural sciences, American University of Ras Al Khaimah, PO Box 10021, Ras Al Khaimah, UAE. E-mail address:
{\tt{ali.maalaoui@aurak.ac.ae}}}

\maketitle

\begin{abstract}
In this paper we provide a method to study critical points of strongly indefinite functionals on vector bundles. We focus mainly on energy functionals coupled with a fermionic part, that is with a Dirac-type operator. We consider the cases of the perturbed Dirac-Geodesics problem and the Yang-Mills-Dirac type equation in dimension two.
\end{abstract}

\section{Introduction and Main Results}

In most of the mathematical physics models involving super-symmetry, the total energy functional involves two parts, a Bosonic classical part and a fermionic part involving a coupling with the Dirac operator. For instance, we can see the Dirac-Harmonic Maps \cite{CJLW,C,Del} and in particular the Dirac-geodesics problem \cite{I,IM}, the Dirac-Einstein functional in full generality, see \cite{Fin, Kim} or under conformal restriction \cite{MV3}, The Yang-Mills-Dirac equation \cite{Li,Ot,Par}, The super-Liouville equation \cite{J}.
The main difficulty in these problems is the fact that the energy functional is strongly indefinite and depending on the dimension, it can be critical. We will focus on the earlier aspect of the problems, that is the strongly indefinite aspect of these energy functionals. This issue comes from the fact that the Dirac operator has infinitely many positive and negative eigenvalues.
There was an extensive work dealing with such problems, involving different  methods. For instance we can cite \cite{M,MV1,MV2,IsobeII} for methods involving a Floer type homology, or \cite{Pan,Sul,Po} for methods involving the generalized Nehari manifold. In this paper, we will rely mainly on the last type of methods. In certain cases, particularly the ones that we will consider, one cannot define the full Nehari manifold as in the classical case, so we will consider here the "half" generalized Nehari manifold to handle the spinoial part of the functional. \\

As a first application of our method, we consider the Dirac-geodesic problem. This problem is the one dimensional version of the perturbed Dirac-Harmonic maps problem, which appears in the non-linear super-symmetric Sigma model (see \cite{Del}). That is we consider the functional
$$E(\phi,\psi)=\frac{1}{2}\int_{S^{1}}\Big|\frac{d\phi}{ds}\Big|^{2}\,ds+\frac{1}{2}\int_{S^{1}}\langle\psi,\D_{\phi}\psi\rangle\,ds-\int_{S^{1}}K(s,\phi(s),\psi(s))\,ds.$$
We show that
\begin{theorem}
Given a compact closed Riemannian manifold $N$, under the assumptions $(H1)-(H4)$, the Dirac-Geodesic problem has infinitely many non-trivial solutions on each homotopy class $[\alpha]\in \pi_{1}(N)$
\end{theorem}

Next, we consider another super-symmetric model, namely the Yang-Mills-Dirac problem in dimension two. Indeed, given a Spin Riemann surface $(M,g,\Sigma)$ and a compact Lie group defining a principal bundle $\pi:P\to M$, we consider the functional 
$$YMD(A,\psi)=\int_{M}|F_{A}|^{2}dv+\frac{1}{2}\int_{M}\langle D_{A}\psi,\psi\rangle dv -\int_{M}K(\psi)dv.$$
Then we have
\begin{theorem}
If $K$ satisfies $(HK)$, then the functional $YMD$ has infinitely many non-gauge-equivalent, non-trivial critical points. 
\end{theorem}

\section{General Setting}

We consider a functional $E:H\to \R$ such that $\pi:H\to M$ is a vector bundle with fibers modeled on the Hilbert space space $H_{u}=V$. From now on, we will drop the subscript $u$ for the fiber unless it is needed. We assume that $$E(u,v)=E_{1}(u)+E_{2}(u,v),$$
where $$E_{2}(u,v)=\langle L v,v\rangle -b(u,v).$$
Here, the $E_{1}$ represents the bosonic part and $E_{2}$ will be the coupled fermionic part. We will assume that the fermionic part $E_{2}(u,v)$ takes the form $$E_{2}(u,v)=\langle L_{u}v,v\rangle -b(u,v),$$ but since $H$ is a second countable infinite dimensional Hilbert manifold, by theorem of Eells and Elworthy (1970), it can be embedded as an open set of a Hilbert space $\tilde{N}\times V$. Thus, we can assume that 
$$\iota(L_{u})=L+g(u,\cdot),$$
where $\iota$ is the map induced by the embedding $H\subset \tilde{N}\times V$ and $L:V\to V$. So from now on, we will identify these two operators and we will absorb the $g$ part in the $b$ functional.\\

We assume that the Hilbert space $V$ is embedded in a dense and compact way in a Hilbert space $(W,|\cdot|)$ so that the operator
$$L: V\longrightarrow W$$
is invertible and self-adjoint. Hence $L$ will have a basis of eigenfunctions $\{\varphi_{i}\}_{i\in \mathbb{Z}}$
$$L(\varphi_i)=\lambda_i \varphi_i$$
with the convention that if $\lambda_i>0$ then $i>0$. This allows us to define the unbounded operator $|L|^{\frac{1}{2}}$ in the following way: if
$$v=\sum_{i\in \mathbb{Z}}a_{i}\varphi_{i}$$
then
$$L(v)=\sum_{i\in \mathbb{Z}}\lambda_{i}a_{i}\varphi_{i}$$
and therefore
$$|L|^{\frac{1}{2}}v=\sum _{i\in \mathbb{Z}}|\lambda_{i}|^{\frac{1}{2}}a_{i}\varphi_{i}.$$
Now if we denote $\langle \cdot,\cdot \rangle$ the inner product in $W$, we define then the inner product of $V$ as follows
$$\langle v_{1},v_{2} \rangle_{V} = \langle |L|^{\frac{1}{2}}v_{1},|L|^{\frac{1}{2}}v_{2} \rangle$$
We obtain the decomposition
$$V=V^{+}\oplus V^{-},$$
where
$$V^{-}=\overline{span\{\varphi_{i},i<0 \}}^{V}, \qquad V^{+}=\overline{span\{\varphi_{i},i>0 \}}^{V}$$
We will write
$$v=v^{+}+v^{-}, \qquad \forall \; v\in V$$
according to the previous splitting also we will write $P^{+}:V\to V^{+}$ and $P^{-}:V\to V^{-}$ the orthogonal projectors on their respective spaces. We explicitly note that
$$L(v^{+} + v^{-}) = |L|(v^{+} - v^{-}).$$
Therefore we will write $\langle |L|v,v \rangle$ in place of $\|v^+\|_V^2 + \|v^-\|_V^2$. It is important to point out here that this way we can construct a two vector bundles $H^{+}$ and $H^{-}$ on $M$ since we can do this splitting at every point of $u\in M$ and the splitting varies smoothly and they are defined as
$$H^{+}=\cup_{u\in M}V_{u}^{+} \text{ and } H^{-}=\cup_{u\in M} V_{u}^{-}.$$
The functional $b$ will be assumed to be compact and $C^{2}$ and such that $$\nabla_{v}b(u,v)=f(u,|v|)v$$ and

\begin{itemize}
\item[i)] $\langle \nabla_{v}b,v\rangle -2b(u,v)>C_{1}(\|u\|)|v|^{p+1}$
\item[ii)] $|\nabla^{2}_{vv}b(u,v)|\leq C_{2}(\|u\|)|v|^{p-1}$
\item[iii)] $f(u,\cdot)$ is increasing and $f(u,0)=0$
\item[iv)] $b(u,v)\geq b(u,0)=0$ and $\frac{b(u,sv)}{s^{2}}\to \infty$ for all $v\not=0$.
\end{itemize}
where we mean by $C_{i}(\|u\|)$ that the constant depends continuously on the magnitude of $u$ and $p>1$ and $\|u\|$ is to be understoud as the distance in $M$ with respect to a fixed reference point $u_{0}$.\\
We recall that a $C^{1}$ functional $F:X\to \R $, where $X$ is a Banach-Finsler manifold, is said to satisfy the Palais-Smale condition (PS), if for every sequence $(x_{n})\subset X$ such that $F(x_{n})\to c$ and $\partial F(x_{n})\to 0$ (such sequence will be called a (PS) sequence), then we can extract a convergent subsequence from $(x_{n})$. This condition is fundamental in the study of variational problem since it is the main ingredient for the classical deformation Lemma.\\
 
We define the generalized Nehari manifold by
$$\mathcal{N}=\Big \{(u,v)\in H\setminus H^{-}, \langle Lv,v\rangle =\langle \nabla_{v}b,v\rangle; P^{-}(Lv-\nabla_{v}b)=0\Big\}.$$
Then one have

\begin{lemma}
The set $\mathcal{N}$ is a manifold.
\end{lemma}
{\it Proof:}
We consider the map $G:H\setminus H^{-}\to \R\times H^{-}$ defined by
$$G(u,v)=\left[\begin{array}{ll}
\langle Lv,v\rangle -\langle \nabla_{v}b,v\rangle\\
 (u,P^{-}(Lv-\nabla_{v}b))
\end{array}
\right].$$
Then clearly $\mathcal{N}=G^{-1}(0)$ hence, if we can show that $d_{v}G(u,v)$ is onto for every $(u,v)\in \mathcal{N}$, we can deduce that the last set is a manifold, since the $u$ component is untouched. For this matter, we restrict our variations first to the $v$ component. So that
$$d_{v}G(u,v)[h_{1},h_{2}]= \left[\begin{array}{ll}
2\langle Lv,h_{1}\rangle -\langle \nabla_{vv}^{2}b h_{1},v\rangle-\langle \nabla_{v}b,h_{1}\rangle\\
 P^{-}(Lh_{2}-\nabla_{vv}^{2}bh_{2})
\end{array}
\right].$$
Hence, if $h_{1}=tv$ and $h_{2}\in V^{-}$, we have that
$$d_{v}G(u,v)[h_{1},h_{2}]= \left[\begin{array}{ll}
t(2\langle Lv,v\rangle -\langle \nabla_{vv}^{2}b v,v\rangle-\langle \nabla_{v}b,v\rangle\\
P^{-}(Lh_{2}-\nabla_{vv}^{2}bh_{2})
\end{array}
\right].$$
But since $(u,v)\in \mathcal{N}$, we have that
$$t(2\langle Lv,v\rangle -\langle \nabla_{vv}^{2}b v,v\rangle-\langle \nabla_{v}b,v\rangle=t( \langle \nabla_{v}b,v\rangle-\langle \nabla_{vv}^{2}b v,v\rangle).$$
Hence, from $iii)$ we have that
$$\langle \nabla_{v}b,v\rangle-\langle \nabla_{vv}^{2}b v,v\rangle<0$$ 
and on $V^{-}$, we have that
$$\langle P^{-}(Lh_{2}-\nabla_{vv}^{2}bh_{2}),h_{2}\rangle =-\|h_{2}\|^{2}-\langle \nabla_{vv}^{2}bh_{2},h_{2}\rangle,$$
which is a negative defined operator, hence invertible. Therefore, we have that $d_{v}G(u,v):\R v\oplus V^{-}\to \R\times V^{-}$ is onto for all $(u,v)\in \mathcal{N}$.  

\hfill$\Box$

We define the set $F_{u}(v)=\R^{+} v\oplus V^{-}$. 

\begin{proposition}
For every $(u,v)\in H\setminus H^{-}$ there exists a unique $v_{0}\in F_{u}(v)$ such that $(u,v_{0})\in \mathcal{N}$.
\end{proposition}
{\it Proof:}
First we show that $E_{2}$ has a maximum on $F_{u}(v)$. So we start by claiming that there exists $R>0$ such that $E_{2}(u,w)\leq 0$ when $\|w\|_{V}>R$.  So we reason by contradiction assuming that there exists a sequence $w_{n}\in F_{u}(v)$ such that $\|w_{n}\|_{V}\to \infty$ and $E_{2}(u,w_{n})>0$. Without loss of generality we can assume that $v=v^{+}\in V^{+}$ and $\|v\|_{V}=1$ since $F_{u}(tv)=F_{u}(v^{+})=F_{u}(v)$. Then we can write $w_{n}=t_{n}v+\varphi_{n}$ and $\|w_{n}\|_{V}^{2}=t_{n}^{2}+\|\varphi_{n}\|_{V}^{2}$. We set $h_{n}=\frac{w_{n}}{\|w_{n}\|_{V}}=s_{n}v+\psi_{n}$. 

Notice that since $s_{n}$ and $\|\psi_{n}\|_{V}$ are bounded, we have that up to a subsequence, $s_{n}\to s_{0}$ and $\psi_{n}\rightharpoonup \psi_{0}$. Therefore,
$$\frac{E_{2}(u,w_{n})}{\|w_{n}\|_{V}^{2}}=\frac{1}{2}(s_{n}^{2}-\|\psi_{n}\|_{V}^{2})-\frac{b(u,w_{n})}{\|w_{n}\|_{V}^{2}}.$$
Thus
$$\limsup \frac{E_{2}(u,w_{n})}{\|w_{n}\|_{V}^{2}}\leq \frac{1}{2}(s_{0}^{2}-\|\psi_{0}\|_{V}^{2})-\liminf \frac{b(u,w_{n})}{\|w_{n}\|_{V}^{2}},$$
but
$$\frac{b(u,w_{n})}{\|w_{n}\|_{V}^{2}}=\frac{b(u,\|w_{n}\|h_{n})}{\|w_{n}\|^{2}_{V}}.$$ 
Therefore, if $h_{n}\rightharpoonup 0$ then $s_{0}\to 0$ and $\psi_{n}\rightharpoonup 0$ thus
$$\limsup \frac{E_{2}(u,w_{n})}{\|w_{n}\|_{V}^{2}}\leq 0,$$
leading to a contradiction. On the other hand, if $h_{n}\rightharpoonup h_{0}\not=0$ hence by $(iv)$, we have that
$$\limsup \frac{E_{2}(u,w_{n})}{\|w_{n}\|_{V}^{2}}\to -\infty.$$
Which leads again to a contradiction. Therefore we can set $$\beta=\sup_{F_{u}(v)}E_{2}(u,\cdot).$$
We claim that $\beta>0$. Indeed, using a Taylor expansion around zero we have, 
$$E_{2}(u,tv)=\frac{1}{2}t^{2}-b(u,tv)=\frac{1}{2}t^{2}-o(t^{2}).$$
Hence we see that for $t>0$ small enough, we have that $E_{2}(u,tv)>0$ and thus $\beta>0$.

So we consider now a maximizing sequence $w_{n}=t_{n}v+\varphi_{n}\in F_{u}(v)$. Clearly, $\|w_{n}\|_{V}$ us bounded, So we can extract again a subsequence, such that $w_{n}\rightharpoonup w_{0}=t_{0}v+\varphi_{0}$. But 
\begin{align}
-\beta&=-\limsup E_{2}(u,w_{n})\notag\\
&\geq \frac{1}{2}(\|\varphi_{0}\|_{V}^{2}-t_{0}^{2})+\liminf b(u,w_{n}).\notag
\end{align}
So by compactness of $b$, we have that
$$-\beta \geq -E_{2}(u,w_{0}).$$
Whence, $E_{2}(u.w_{0})=\beta$ and we do indeed have a maximize and we need to show now the uniqueness of the maximizer. So let us take $(u,v)\in \mathcal{N}$
We want to show that $E_{2}(u,tv+w)<E_{2}(u,v)$ unless $t=1$ and $w=0$. In fact, one has
$$E_{2}(u,tv+w)=\frac{1}{2}(t^{2}\langle Lv,v\rangle-\|w\|_{V}^{2}+2t\langle Lv,w\rangle)-b(u,tv+w).$$
But since $(u,v)\in\mathcal{N}$ we have that
$$E_{2}(u,tv+w)=\langle \nabla_{v}b(u,v),\frac{1}{2}t^{2}v+tw \rangle-\frac{1}{2}\|w\|_{V}^{2}-b(u,tv+w).$$
Hence,
$$E_{2}(u,tv+w)-E_{2}(u,v)=\langle \nabla_{v}b(u,v),\frac{1}{2}(t^{2}-1)v+tw \rangle +b(u,v)-b(u,tv+w)-\frac{1}{2}\|w\|_{V}^{2}.$$
In particular if $h(t)=\langle \nabla_{v}b(u,v),\frac{1}{2}(t^{2}-1)v+tw \rangle +b(u,v)-b(u,tv+w)$ is negative then we have the desired result. This last claim of negativity follows exactly from the procedure in \cite{Sul} and \cite{Po}.

\hfill$\Box$.

For $(u,v)\in H\setminus H^{-}$, we will denote by $g_{u}(v)=s_{u}(v)v+\varphi_{u}(v)$ the map such that $(u,g_{u}(v))\in \mathcal{N}$. Notice that since $\mathcal{N}$ is a manifold is equivalent to the smoothness of the map $g$. We define thus the functional $$\tilde{E}(u,v)=E(u,g_{u}(v)).$$

\begin{lemma}
If $E_{1}$ is coercive, then any  Palais-Smale sequence $(u_{n},v_{n})$ of $E_{|\mathcal{N}}$,  is a Palais-Smale sequence of $E$.
\end{lemma}

{\it Proof}
Notice first that $$\tilde{E}(u,v)=E(u,g_{u}(v))=E_{1}(u)+\frac{1}{2}\langle \nabla_{v}b(u,g_{u}(v))\rangle -b(u,g_{u}(v)).$$
Therefore, from $i)$, it is bounded from below, hence if $(z_{n})$ is a (PS) sequence for $\tilde{E}$, then $\|u_{n}\|$ is bounded and so is $|g_{u_{n}}(v_{n})|$.

Now, we have that
\begin{align}
\partial_{v}\tilde{E}(u,v)[h]&=\nabla_{v} E_{2}(u,g_{u}(v))[\partial_{v}g_{u}(v)[h]]\notag \\
&=\partial_{v} E_{2}(u,g_{u}(v))[\partial_{v}t_{u}(v)[h]v+t_{u}(v)h+\partial_{v}\varphi [h]]\notag\\
&=t_{u}(v)\partial_{v} E_{2}(u,g_{u}(v)[h]\notag 
\end{align}
and $$t_{u}(v)=\frac{\|g_{u}(v)^{+}\|_{V}}{\|v^{+}\|_{V}}.$$
On the other hand
$$\partial_{u}\tilde{E}(u,v)[h]=\partial_{u}E_{1}(u)[h]+\partial_{u} E_{2}(u,g_{u}(v))[h]+\partial_{v} E_{2}(u,g_{u}(v))[\partial_{u}g_{u}(v)[h]]\rangle$$
$$=\partial_{u}E_{1}(u)[h]+\partial_{u} E_{2}(u,g_{u}(v))[h].$$
Hence, if $(u_{n},v_{n})\in \mathcal{N}$ is a (PS) sequence of $E_{|\mathcal{N}}$, as long as $\|v^{+}\|_{V}$ is bounded away from zero, we do have indeed a (PS) sequence for $E$.
notice that if $(u,v)\in \mathcal{N}$ then we have that
\begin{equation}\label{in1}
-\|v^{-}\|^{2}_{V}=\langle \nabla_{v}b(u,v),v^{-}\rangle.
\end{equation}
Thus,
\begin{equation}\label{in2}
\|v^{-}\|_{V}\leq f(u,|v|)|v|,
\end{equation}
also
\begin{equation}\label{in3}
\|v\|_{V}^{2}=\langle Lv,v\rangle +2\|v^{-}\|_{V}^{2}\leq f(u,|v|)(1+f(u,|v|))|v|^{2}.
\end{equation}
Whence
$$1\leq C(f(u,|v|)(1+f(u,|v|))).$$
Letting $v\to 0$ we find a contradiction. hence $\|v\|_{V}>\delta$.
On the other hand, we have that
$$\|v^{+}\|_{V}^{2}-\|v^{-}\|_{V}^{2}=\langle \nabla_{v}b(u,v),v\rangle\geq 0.$$
Thus, $\|v^{-}\|_{V}^{2}\leq \|v^{+}\|_{V}^{2}$, and therefore
$$\delta<\|v\|_{V}^{2}=\|v^{+}\|_{V}^{2}+\|v^{-}\|_{V}^{2}\leq 2\|v^{+}\|_{V}^{2}.$$
 Hence any (PS) sequence of $\tilde{E}$ is a (PS) sequence of $E$.
\hfill$\Box$

\begin{lemma}
If $E_{1}$ is coercive and weakly lower semi-continuous then $E$ has at least one critical point.
\end{lemma}
{\it Proof}
Let us consider a minimizing sequence of $E_{|\mathcal{N}}$, then by coercivity of $E_{1}$ we have that $\|u_{n}\|$ is bounded and hence it converges weakly to $u_{\infty}$. This also implies the boundedness of $|v_{n}|^{p+1}$ and using inequalities $(\ref{in2})$ and $(\ref{in3})$, we have the boundedness of $\|v_{n}\|_{V}$. Thus, there exist a weakly convergent subsequence that converges to weakly $v_{\infty}$ in $H$ and strongly in $W$. Now if $(u_{\infty},v_{\infty})\in \mathcal{N}$ then we do have a minimizer, which will be a critical point of $E$. 

Since 
$$\langle Lv_{n}-\nabla_{v}b(u_{n},v_{n}),\varphi\rangle =0, \text{ for all } \varphi\in V^{-},$$
by passing to the limit, we have that
$$P^{-}( Lv_{\infty}-\nabla_{v}b(u_{\infty},v_{\infty}))=0.$$
Moreover, since $z_{n}$ is a (PS) sequence for $E$, we have in particular that
$$Lv_{n}-\nabla_{v}b(u_{n},v_{n})=o(1).$$
Testing against $v_{\infty}$ we see that 
$$\langle Lv_{n}-\nabla_{v}b(u_{n},v_{n}),v_{\infty}\rangle=o(1)$$
using the weak convergence and passing to the limit, we see that $(u_{\infty},v_{\infty})$ is indeed in $\mathcal{N}$, moreover, we do have the strong convergence of $v_{n}\to v_{\infty}$ and hence
$$E_{|\mathcal{N}}(u_{\infty},v_{\infty})\leq \liminf E(u_{n},v_{n}).$$
So we have indeed one non-trivial critical point.
\hfill$\Box$

\section{Coupling with the Dirac Operator}
Given a Riemannian spin manifold $M$, we let $\Sigma M$ denote the canonical spinor bundle associated to $M$, see \cite{F}, whose sections are simply called spinors on $M$. This bundle is endowed with a natural Clifford multiplication $c$, a hermitian metric and a natural metric connection $\nabla^\Sigma$. The Dirac operator $D_g$ acts on spinors 
$$D_g:C^\infty(\Sigma M)\longrightarrow C^\infty(\Sigma M)$$ defined as the composition $c \circ \nabla^\Sigma$, where $c$ is the Clifford multiplication, in the following way: if $(e_{1},\cdots,e_{n})$ is an orthonormal local frame of $TM$, then $$D_{g}\psi=\sum_{i=1}^{n}e_{i}\cdot \nabla_{e_{i}}^{\Sigma}\psi.$$

The functional space that we are going to define is the Sobolev space $H^{\frac{1}{2}}(\Sigma M)$. First we recall that the Dirac operator $D_g$ on a compact manifold is essentially self-adjoint in $L^2(\Sigma M)$ and has compact resolvent and there exists a complete $L^2$-orthonormal basis of eigenspinors $\{\psi_i\}_{i\in\mathbb{Z}}$ of the Dirac operator
$$D_g\psi_i=\lambda_i \psi_i ,$$
and the eigenvalues $\{\lambda_i\}_{i\in\mathbb{Z}}$ are unbounded, that is $|\lambda_i|\rightarrow\infty$, as $i\rightarrow\infty$.
Now if $\psi\in L^{2}(\Sigma M)$, it has a representation in this basis , namely: 
$$\psi=\sum_{i\in \mathbb{Z}}a_{i}\psi_{i}.$$
Let us define the unbounded operator $|D_g|^{s}: L^{2}(\Sigma M)\rightarrow L^{2}(\Sigma M)$ by 
$$|D_g|^{s}(\psi)=\sum_{i\in \mathbb{Z}} a_{i}|\lambda_{i}|^{2s}\psi_{i}.$$
We denote by $H^s(\Sigma M)$ the domain of $|D_g|^{s}$, namely $\psi\in H^s(\Sigma M)$ if and only if 
$$\sum_{i\in \mathbb{Z}} a_{i}^2|\lambda_{i}|^{2s}<+\infty .$$
$H^s(\Sigma M)$ coincides with the usual Sobolev space $W^{s,2}(\Sigma M)$ and for $s <0$, $H^s(\Sigma M)$ is defined as the dual of $H^{-s}(\Sigma M)$. For s >0, we can define the inner product 
$$\langle u,v\rangle_{s}=\langle|D_g|^{s}u,|D_g|^{s}v\rangle_{L^{2}},$$ 
which induces an equivalent norm in $H^{s}(\Sigma M)$; we will take 
$$\|u\|^{2}=\langle u,u\rangle_{\frac{1}{2}}$$
as our standard norm for the space $H^{\frac{1}{2}}(\Sigma M)$. Even in this case, the Sobolev embedding theorems say that there is a continuous embedding for $dim(M)=n>1$
$$H^{s}(\Sigma M) \hookrightarrow L^p(\Sigma M), \quad 1\leq p \leq \frac{2n}{n-1},$$
which is compact if $1\leq p <\frac{2n}{n-1}$. For $n=1$ we have that the embedding is compact in all $L^{p}$ for $1\leq p<\infty$.

Finally, we will decompose $H^{\frac{1}{2}}(\Sigma M)$ in a natural way. Let us consider the $L^2$-orthonormal basis of eigenspinors $\{\psi_i\}_{i\in\mathbb{Z}}$: we denote by $\psi_i^-$ the eigenspinors with negative eigenvalue, $\psi_i^+$ the eigenspinors with positive eigenvalue and $\psi_i^0$ the eigenspinors with zero eigenvalue; we also recall that the dimension of the kernel of $D_g$ is finite dimensional. Now we set:
$$H^{\frac{1}{2},-}:=\overline{\text{span}\{\psi_i^-\}_{i\in\mathbb{Z}}},\quad
H^{\frac{1}{2},0}:=\overline{\text{span}\{\psi_i^0\}_{i\in\mathbb{Z}}}, \quad
H^{\frac{1}{2},+}:=\overline{\text{span}\{\psi_i^+\}_{i\in\mathbb{Z}}},$$
where the closure is taken with respect to the $H^{\frac{1}{2}}$-topology. Therefore we have the orthogonal decomposition  $H^{\frac{1}{2}}(\Sigma M)$ as:
$$H^{\frac{1}{2}}(\Sigma M)=H^{\frac{1}{2},-}\oplus H^{\frac{1}{2},0}\oplus H^{\frac{1}{2},+}.$$
We will let $P^{+}$, $P^{0}$ and $P^{-}$ be the projectors on $H^{\frac{1}{2},+}$, $H^{\frac{1}{2},0}$ and $H^{\frac{1}{2},-}$.

\subsection{The Dirac-Geodesic Problem}
In this section we will adapt the method stated above to find solutions to the Dirac-Geodesic problem studied in \cite{I,IM}. In fact the proof that we provide here is shorter and much simpler than the one in $\cite{I}$ even though, we deal with a certain class of non-linearities. But we believe that this method can be extended even more to incorporate the cases in $\cite{I}$.

Let  $N$ be a compact Riemannian manifold. We define the configuration space $\mathcal{F}^{1,1/2}(S^{1},N)$ as
$$\mathcal{F}^{1,1/2}(S^{1},N)=\Big \{(\phi,\psi):\phi\in H^{1}(S^{1},N),\,\psi\in H^{1/2}(S^{1},\Sigma S^{1}\otimes\phi^{\ast}TN)\Big\}.$$
This space is disconnected and the connected components are coming for the homotopy classes of the loops $\phi:S^{1}\to N$. Hence, we will restrict the study to each homotopy class $[\alpha]\in \pi_{1}(N)$. Again here, as we saw above, the space $H^{\frac{1}{2}}(S^{1},\Sigma S^{1}\otimes \phi^{*}TN)$ splits into two parts $H^{\frac{1}{2},-}_{\phi}$ and $H^{\frac{1}{2},+}_{\phi}$. We will write then $P_{\phi}^{\pm}$ the projector on $H^{\frac{1}{2},\pm}_{\phi}$  
The operator $D_{\phi}$ is constructed in the following way: First we consider the connection induced by the metrics on $\Sigma S^{1}$ and $\phi^{*}TN$. Then using this connection, we define the Dirac operator by composing with the Clifford multiplication. Indeed, if $D_{0}=i\frac{\partial}{\partial s}$ the untwisted Dirac operator on $\Sigma S^{1}$, and $\psi(s)=\psi^{i}\otimes \frac{\partial}{\partial y_{i}}(\phi(s))$ then the Dirac operator can be expressed locally by
\begin{equation}\label{Dirac}
D_{\phi}\psi=D_{0}\psi^{i}\otimes \frac{\partial}{\partial y_{i}}(\phi(s))+\Gamma^{i}_{jk}(\phi(s)) \frac{\partial \phi }{\partial s}^{j}\cdot \psi^{k}(\phi)\otimes  \frac{\partial}{\partial y_{i}}(\phi(s)),
\end{equation}
where $\Gamma^{i}_{ik}$ are the Christoffel symbols of $N$.
We consider the perturbed Dirac-geodesic action $E$ defined by
$$E(\phi,\psi)=\frac{1}{2}\int_{S^{1}}\Big|\frac{d\phi}{ds}\Big|^{2}\,ds+\frac{1}{2}\int_{S^{1}}\langle\psi,\D_{\phi}\psi\rangle\,ds-\int_{S^{1}}K(s,\phi(s),\psi(s))\,ds,$$
where $K:S^{1}\times\Sigma S^{1}\otimes TN\to\R$ is a smooth function (we write $K=K(s,\phi,\psi)$), where $s\in S^{1}$ and $(\phi,\psi)\in\Sigma S^{1}\otimes TN$, i.e., $\phi\in N$ is a base point and $\psi\in\Sigma S^{1}\otimes T_{\phi}N$ is a point on the fiber over $\phi\in N$). We assume that there exists $p>2$ such that $K$ satisfies
\begin{itemize}
\item[H1)] $|d^{2}_{\psi\psi}K(s,\phi,\psi)|\le C_{1}(1+|\psi|^{p-1}),$

\item[H2)] $C_{2}|\psi|^{p+1}+2K(s,\phi,\psi)\le\langle\nabla_{\psi}K(s,\phi,\psi),\psi\rangle,$

\item [H3)] $\nabla_{\psi}K(s,\phi,\psi)=f(s,\phi,|\psi|)\psi$ and $f$ is increasing with $f(s,\phi,0)=0$,

\item [H4)] $K(s,\phi,\psi)\geq K(s,\phi,0)=0$ and $\frac{K(s,\phi,\lambda \psi)}{\lambda^{2}}\to \infty$ as $\lambda\to \infty$ and $\psi\not=0$.
\end{itemize}

In $\cite{I}$, Isobe proved that
\begin{proposition}[\cite{I}] For $(\phi,\psi)\in\mathcal{F}^{1,1/2}(S^{1},N)$, we have the following:
\begin{equation}
\left\{\begin{array}{ll}
\nabla_{\phi}E(\phi,\psi)=(-\Delta+1)^{-1}\Big(-\nabla_{s}\partial_{s}\phi
+\frac{1}{2}R(\phi)\langle\psi,\partial_{s}\phi\cdot\psi\rangle-\nabla_{\phi}K(s,\phi,\psi)\Big),\\
\nabla_{\psi}E(\phi,\psi)=(1+|\D|)^{-1}(\D_{\phi}\psi-\nabla_{\psi}K(s,\phi,\psi)),
\end{array}
\right.
\end{equation}
where 
$$R(\phi)\langle\psi,\partial_{s}\phi\cdot\psi\rangle=\Big\langle\psi,\partial_{s}\cdot\psi^{i}\otimes\frac{\partial}{\partial y^{j}}(\phi)\Big\rangle\partial_{s}\phi^{l}R^{j}_{iml}(\phi)g^{ms}(\phi)\frac{\partial}{\partial y^{s}}(\phi).$$
\end{proposition}
See~\cite{I} for the details of the derivation of the above formula. So we propose in this case to find solutions to the system
\begin{equation}\label{DE}
\left\{\begin{array}{ll}
-\nabla_{s}\partial_{s}\phi+\frac{1}{2}R(\phi)\langle\psi,\partial_{s}\phi\cdot\psi\rangle=\nabla_{\phi}K(s,\phi,\psi)\\
\D_{\phi}\psi=\nabla_{\psi}K(s,\phi,\psi).
\end{array}
\right.
\end{equation}
Notice that this system has already trivial solutions if we take $\psi=0$ and $\phi$ a geodesic on $N$. Similarly to what we have defined above, we consider the generalized Nehari manifold $\mathcal{N}$ defined by
$$\mathcal{N}=\Big \{(\phi,\psi)\in \mathcal{F}^{1,\frac{1}{2}}(S^{1},N); \int_{S^{1}}\langle D_{\phi}\psi,\psi\rangle =\int_{S^{1}}\langle \nabla_{\psi}K,\psi\rangle; P^{-}_{\phi}(D_{\phi}\psi-\nabla_{\psi}K)=0\Big \}.$$
As we saw above, we can show that $\mathcal{N}$ is indeed a manifold and any (PS) sequence for $E_{|\mathcal{N}}$ is also a (PS) sequence for $E$. It is important to notice here that there is a small but relevant difference, form the case above. In fact, in the above case, the operator $L$ is independent of $u$, but in this case we can take it to be dependent on $\phi$. It appears to be more convenient to do it that way but it does not change any thing to the proof.

First, notice that
$$E_{|\mathcal{N}}(\phi,\psi)=\frac{1}{2}\int_{S^{1}}|\dot{\phi}|^{2}+\frac{1}{2}\int_{S^{1}}\langle \nabla_{\psi}K,\psi\rangle -2K(\phi,\psi),$$
which is bounded from below. 
\begin{lemma}
Let $(z_{n})$ be a Palais-Smale sequence for $E_{|\mathcal{N}}$ then there exists $\delta>0$ such that $\|\psi\|_{H^{\frac{1}{2}}}\geq \delta$.
\end{lemma}
{\it Proof:}
First, notice that since $\phi_{n}$ is bounded in $H^{1}(S^{1},N)$, in particular $\dot{\phi}_{n}$ is bounded in $L^{2}$, we have that the norms defined on the bundle above $\phi_{n}$ is equivalent to the standard one. In fact this follows from the expression $(\ref{Dirac})$ that
$$D_{\phi}=D_{0}+A(\dot{\phi}),$$
where $A(\dot{\phi})$ is linear in $\dot{\phi}$. Now, we have that
$$\|\psi^{+}\|_{\frac{1}{2}}^{2}-\|\psi^{-}\|_{\frac{1}{2}}^{2}\leq C\|\psi\|_{p+1}^{p+1}.$$
On the other hand, we have that
$$-\|\psi^{-}\|_{\frac{1}{2}}^{2}=\int_{M}\langle \nabla_{\psi}H,\psi^{-}\rangle.$$
Hence,
$$\|\psi^{-}\|_{\frac{1}{2}}^{2}\leq C\|\psi\|_{p+1}^{p+1},$$
therefore 
$$\|\psi\|_{\frac{1}{2}}^{2}\leq C\|\psi\|_{p+1}^{p+1}.$$
But from the classical Sobolev embedding, we have that
$$\|\psi\|_{p+1}^{2}\leq C\|\psi\|_{\frac{1}{2}}^{2}\leq C_{1}\|\psi\|^{p+1}_{p+1}.$$
Since, $p>1$,  $\|\psi\|_{\frac{1}{2}}$ cannot converge to zero.
\hfill$\Box$

Now we consider a minimizing sequence $(z_{n})$ of $E_{|\mathcal{N}}$ it follows from Ekeland's variational principle \cite{Ek}, that it is a (PS) sequence for $E$ and since in this case $E$ satisfies the (PS) condition, one has a minimizer, in each homotopy class $[\alpha]\in \pi_{1}(M)$.
In fact, in this case, $E$ satisfies the (PS) condition, (see \cite{I}), then so does $E_{|\mathcal{N}}$. We have then the following result
\begin{theorem}
If we assume moreover that $K$ is even in $\psi$, then we have infinitely many solutions to $(\ref{DE})$.
\end{theorem}
{\it Proof:}

Notice that in this case $\mathcal{N}$ is invariant under the action of $\Z_{2}$ on the $\psi$ component, we consider then $\mathcal{B}_{k}$ the collection of sets $B\subset \mathcal{N}$ such that $-B=B$ and $\gamma(B)\geq k$ where $\gamma$ is the Krasnoselskii genus, also we consider the sequence of numbers $c_{k}$ defined by
$$c_{k}= \inf_{B\in \mathcal{B}_{k}}\max_{B}E.$$
Then we already know from classical min-max theory (see \cite{Rab}), that the $c_{k}$ are critical values of $E$, so if we show that $\gamma(\mathcal{N})=\infty$, we do have indeed infinitely many solutions. So we fix $\phi\in H^{1}(S^{1},N)$ and we consider the map $\tilde{T}:S^{+}_{\phi}\to \mathcal{N}$ defined by $$\tilde{T}(\psi)=(\phi,T_{\phi}(\psi)),$$
where $S_{\phi}^{+}$ is the unit sphere of $H^{\frac{1}{2},+}_{\phi}$. Then by uniqueness of the maximizer as in Proposition 2.2, we have that 
$$\tilde{T}(-\psi)=(\phi, -T_{\phi}(\psi)).$$
Since $\gamma(S^{+})=\infty$, we have then $\gamma(\mathcal{N})=\infty$ leading to the desired result.
\hfill$\Box$

\section{The Yang-Mills-Dirac Problem}
In this section we consider a Riemann surface $(M,g)$ and a compact Lie group $G$ with principal $G$-bundle $P\to M$. If $\sigma:G\to Aut(\mathfrak{g})$ is the adjoint representation of $G$, we define the adjoint vector bundle $Ad(P)=P\times_{\sigma}\mathfrak{g}$. A smooth connection $A$ on $P$ is an equivariant $\mathfrak{g}$-valued $1$-form, with values in the vertical direction, that is $A\in \Omega^{1}(P,\mathfrak{g})$ satisfying for $p\in P$, $v\in T_{p}P$, $h\in G$ and $\xi\in \mathfrak{g}$,
\begin{itemize}
\item $A_{ph}(vh)=h^{-1}A_{p}(v)h$
\item $A_{p}(p\xi)=\xi.$
\end{itemize}
We will set $\mathcal{A}(P)$ the set of smooth connections on $P$. Every connection $A$ on $P$, provides a covariant derivative $\nabla_{A}:C^{\infty}(M,Ad(P))\to C^{\infty}(M,T^{*}M\otimes Ad(P))$ that can be extended to an exterior differential $d_{A}: C^{\infty}(M,\Lambda^{p}T^{*}M\otimes Ad(P))\to C^{\infty}(M,\Lambda^{p+1}T^{*}M\otimes Ad(P))$
Locally, $d_{A}$ can be expressed as $d_{A}=d+\sigma_{*}(A)$. The curvature of a connection is the two form $F_{A}=(d_{A})^{2}$ that we can write as
$$F_{A}=dA+\frac{1}{2}[A,A].$$
One can check that $$F_{A_{1}}=F_{A_{0}}+d_{A_{0}}(A_{1}-A_{0})+\frac{1}{2}[A_{1}-A_{0},A_{1}-A_{0}]$$
and $$d_{A_{1}}-d_{A_{2}}=[A_{1}-A_{2},\cdot],$$
for $A_{1},A_{2}\in \mathcal{A}(P)$. For further details on gauge theory, we refer the reader to \cite{Weh}.

In a similar way as for the Levi-Civita connection, we can extend the connection $\nabla_{A}$ to the bundle $\mathcal{H}=\Sigma M\otimes Ad(P)$ locally by
$$\tilde{\nabla}_{A}(s\otimes v)=\nabla s\otimes v+s\otimes \nabla_{A}v.$$
Hence, one can define the Twisted Dirac operator $D_{A}$ on sections of $\mathcal{H}$ as $D_{A}=c\circ \tilde{\nabla}_{A}$ where $c$ is the Clifford multiplication.
We recall also the Gauge group $\mathcal{G}(P)$, which is the set of equivariant maps $u:P\to G$. The action of the group $\mathcal{G}(P)$ on $\mathcal{A}(P)$ is defined by
$$u^{*}A=u^{-1}Au+u^{-1}du.$$
With this action, we notice that $$F_{u^{*}A}=u^{-1}F_{A}u.$$
Moreover, we can define an action of $\mathcal{G}(P)$ on $\mathcal{H}=\mathbb{S}(M)\otimes Ad(P)$ by
$$u^{*}(s\otimes v)=s\otimes u^{-1}v.$$
With this action, we have that
$$D_{u^{*}A}u^{*}\psi=u^{*}(D_{A}\psi).$$
We can also define the Sobolev Spaces of connections $\mathcal{A}^{k,p}(P)$ as the space of connections in $L^{p}$, with derivatives up to order $k$ in $L^{p}$. In particular, $\mathcal{A}^{1}(P)=\mathcal{A}^{1,2}(P)$ is the substitute of the Sobolev space with Hilbert structure $H^{1}$ In fact, if $(A,B)_{L^{2}}$ defines the $L^{2}$ inner product on $\mathcal{A}^{0,2}(P)$, that is,
$$(A,B)_{L^{2}}=\int_{M}(A,B)dv,$$
then we define the $H^{1}$, inner product with respect to a given connection $A_{0}$ by
$$(A,B)_{A_{0}}=(A,B)_{L^{2}}+(\nabla_{A_{0}}A, \nabla_{A_{0}}B)_{L^{2}}.$$
The associated norm will then be denoted by $\|\cdot\|_{A_{0}}$. The norm on the dual space $\mathcal{A}^{-1}(P)$ will be denoted by $\|\cdot\|_{A_{0}}^{*}$.  Moreover we have the following 
\begin{lemma}[\cite{Wil}]
Let $S$ be a bounded set in $\mathcal{A}^{0,4}(P)$, the set of connections in $L^{4}$. Then if $A_{1}, A_{2}\in S$, there exists $C(S)>0$ depending on the bound of $S$, such that for all $A\in \mathcal{A}^{1}(P)$,
$$C(S)^{-1}\|A\|_{A_{1}}\leq \|A\|_{A_{2}}\leq C(S)\|A\|_{A_{1}}$$
and for all $B\in \mathcal{A}^{-1}(P)$
$$C(S)^{-1}\|B\|_{A_{1}}^{*}\leq \|B\|_{A_{2}}^{*}\leq C(S)\|B\|_{A_{1}}^{*}.$$
\end{lemma}

Also $\mathcal{G}^{2,2}(P)$ the space of maps that are square integrable and with derivatives up to the second order, square integrable (see \cite{Weh} for details). The space $H^{\frac{1}{2}}(\mathcal{H})$ is defined in the usual way as in the introduction of Section 3 with respect to a fixed connection $A_{0}$ and the norm will be denoted by $\|\cdot\|_{\frac{1}{2},A_{0}}$ and the dual norm will be denoted by $\|\cdot\|_{\frac{1}{2},A_{0}}^{*}$. Also one can show easily the following
\begin{lemma}
If $S$ is a bounded set in $\mathcal{A}^{0,q}(P)$ the set of connections in $L^{q}$ for $q>4$, then for $A_{1},A_{2}\in S$, there exists $C(S)>0$ depending on the bound on $S$, such that for $\psi \in H^{\frac{1}{2}}(\mathcal{H})$,
$$C(S)^{-1}\|\psi\|_{\frac{1}{2},A_{1}}\leq \|\psi\|_{\frac{1}{2},A_{2}}\leq C(S)\|\psi\|_{\frac{1}{2},A_{1}}$$
and for all $\varphi\in H^{-\frac{1}{2}}(\mathcal{H})$,
$$C(S)^{-1}\|\varphi\|_{\frac{1}{2},A_{1}}^{*}\leq \|\varphi\|_{\frac{1}{2},A_{2}}^{*}\leq C(S)\|\varphi\|_{\frac{1}{2},A_{1}}^{*}.$$

\end{lemma}

 In these spaces, we can define the functional $YMD:\mathcal{A}^{1}\times H^{\frac{1}{2}}(\mathcal{H})\to \R$ by
$$YMD(A,\psi)=\int_{M}|F_{A}|^{2}dv+\frac{1}{2}\int_{M}\langle D_{A}\psi,\psi\rangle dv -\int_{M}K(\psi)dv$$
Where for sinplicity here we will take 
\begin{equation}\label{hyp}
(HK)\qquad K(x,\psi)=\frac{1}{p+1}b(x)|\psi|^{p+1},
\end{equation} 
with $b$ a smooth strictly positive function on $M$ and $2<p+1<4$.
\begin{proposition}[\cite{Par,Ot}]
The critical points of $YMD$ satisfie the equation
\begin{equation}\label{YD}
\left\{\begin{array}{ll}
\delta_{A}F_{A}=J(\psi,\psi)\\
D_{A}\psi=c(x)|\psi|^{p-1}\psi
\end{array}
\right.
\end{equation}
where $J(\psi,\psi)=-\frac{1}{2}\langle \psi, e^{i}\cdot \sigma(g^{\alpha})\psi\rangle e_{i}\otimes a^{\alpha}$ where $(a_{\alpha})$ is an orthonormal basis of $\mathfrak{g}$ and $(e_{i})$ is a local frame of $TM$ and $\sigma$ is the unitary representation. The operator $\delta_{A}$ is the formal adjoint of $d_{A}$.
\end{proposition}
Again, here we have that the functional $YMD$ is the sum of two functional $YMD(A,\psi)=E_{1}(A)+E_{2}(A,\psi)$, where
$$E_{1}(A)=YM(A)=\int_{M}|F_{A}|^{2}dv$$
and
$$E_{2}(A,\psi)=\frac{1}{2}\int_{M}\langle D_{A}\psi,\psi\rangle dv -\int_{M}K(\psi)dv.$$
We recall that the functional $YM$ was extensively investigated because of its topological and geometrical implications. We refer the reader to \cite{At} for the study of the functional in dimension two and \cite{Don,Fr} for its study in dimension four. Notice that $(\ref{YD})$ has trivial solutions by taking $\psi=0$ and $A$ a Yang-Mills connection, but in this work we are interested in non-trivial solutions, that is $\psi\not=0$.

The space $H^{\frac{1}{2}}(\mathcal{H})$ splits in a natural way with respect to the spectrum of the Dirac operator $D_{A}$ as
$$H^{\frac{1}{2}}(\mathcal{H})=H^{+,\frac{1}{2}}_{A}\oplus H^{0,\frac{1}{2}}_{A}\oplus H^{-,\frac{1}{2}}_{A}.$$
We will also denote by $H^{-,0}_{A}= H^{0,\frac{1}{2}}_{A}\oplus H^{-,\frac{1}{2}}_{A}$ and again $P_{A}^{+}$, $P_{A}^{-}$, $P_{A}^{0}$ the projectors on $H^{+,\frac{1}{2}}_{A}$, $H^{-,\frac{1}{2}}_{A}$ and $H^{0,\frac{1}{2}}_{A}$ respectively. We will also take $P_{A}^{-.0}=P_{A}^{-}+P_{A}^{0}$.\\

Clearly, $YMD$ is invariant under the action of $\mathcal{G}^{2,2}(P)$. We can now define the generalized Nehari manifold by
$$\mathcal{N}=\Big \{(A,\psi)\in H\setminus H^{-.0}; \int_{M}\langle D_{A}\psi,\psi \rangle =\int_{M}\langle K'(\psi),\psi\rangle dv;
P_{A}^{-,0}(D_{A}\psi-K'(\psi))=0\Big\}.$$
Notice that since $H_{u^{*}A}^{\pm,\frac{1}{2}}=u^{-1}H_{A}^{\pm,\frac{1}{2}}$, we deduce that $\mathcal{N}$ is invariant under the action of $\mathcal{G}(P)$. As in the previous sections we define the space $H_{A}(\psi)=\R^{+}\psi\oplus H^{-,0}_{A}$.
\begin{proposition}
Given $A\in \mathcal{A}^{1}(P)$ and $\psi \in H^{\frac{1}{2}}(\mathcal{M})\setminus H_{A}^{-,0}$, then the functional $E_{2}(A,\cdot)_{|H_{A}(\psi)}$ has a unique maximizer $T_{A}(\psi)=t_{A}(\psi)\psi+\varphi_{A}(\psi)$.
\end{proposition}
Notice that from the uniqueness, we have that $$T_{u^{*}A}(u^{*}\psi)=u^{*}T_{A}(\psi).$$
Define then functional $\widetilde{YMD}(A,\psi)=YMD(A,T_{A}(\psi))$.
\begin{lemma}
The Palais-Smale sequences of $YMD_{|\mathcal{N}}$ are Palais-Smale sequences of $YMD$. In particular, the critical points of $YMD_{|\mathcal{N}}$ are also critical points of $YMD$.
\end{lemma}
{\it Proof:}
This follows from the fact that 
\begin{align}
\partial_{\psi} E_{2}(A,T_{A}(\psi))[h]&=(\partial_{\psi}E_{2})(A,T_{A}(\psi))[\partial_{\psi}T_{A}(\psi)[h]]\notag\\
&=(\partial_{\psi}E_{2})(A,T_{A}(\psi))[\partial_{\psi}t_{A}(\psi)[h]\psi+\partial_{\psi}\varphi_{A}(\psi)[h]+t_{A}(\psi)h]\notag\\
&=t_{A}(\psi)(\partial_{\psi}E_{2})(A,T_{A}(\psi))[h],\notag
\end{align}
where $t_{A}(\psi)=\frac{\|T_{A}(\psi)^{+}\|_{\frac{1}{2},A}}{\|\psi^{+}\|}$ and
\begin{align}
\partial_{A}YMD(A,T_{A}(\psi))[h]&=(\partial_{A}E_{1})(A)[h]+(\partial_{A}E_{2})(A,T_{A}(\psi))[h]\notag\\
&\quad +(\partial_{\psi}E_{2})(A,T_{A}(\psi))[\partial_{A}t_{A}(\psi)[h]\psi+\partial_{A}\varphi_{A}(\psi)[h]]\notag\\
&=(\partial_{A}E_{1})(A)[h]+(\partial_{A}E_{2})(A,T_{A}(\psi))[h].\notag
\end{align}
Hence, it is enough to show that there exists $\delta>0$ such that $\|\psi^{+}\|_{\frac{1}{2},A}>\delta$, for all $(A,\psi)\in \mathcal{N}$. Indeed, if $\psi\in S^{+}_{A,r}$, the sphere of radius $r>0$ of $H^{+,\frac{1}{2}}_{A}$, we have that
$$E_{2}(A,\psi)=r^{2}-\int_{M}K(\psi)\geq r^{2}-cr^{p+1},$$
therefore for $r>0$ and small enough we have the existence of $\delta_{1}>0$ such that $E_{2}(A,\psi)>\delta_{1}$. Now, notice that $$E_{2}(A,T_{A}(\psi))=\max_{t>0,\phi \in H_{A}^{-,0}}E_{2}(A,t\varphi^{+}+\phi)\geq \delta_{1}.$$
Thus, there exists $\delta>0$ such that $$\|T_{A}(\psi)^{+}\|_{\frac{1}{2},A}>\delta.$$
\hfill$\Box$

Recall that by taking the space $\mathcal{G}^{2}_{m}$, of gauge transformations fixing the fiber above $m\in M$, then we have that $\mathcal{G}^{2}_{m}$ acts freely on $\mathcal{A}^{1}(P)$ hence, the action is also free on $H$. Thus the space $\overline{H}=H/\mathcal{G}^{2}_{m}$ has the structure of a manifold, moreover the functional $YMD$ descends to the quotient as $\overline{YMD}$ as a well defined functional on $\overline{H}$ and it is $C^{2}$. We can also take the quotient of $\mathcal{N}$ under the action of $\mathcal{G}^{2}_{m}$ that we will denote by $\overline{\mathcal{N}}$. Notice also that $\mathcal{G}^{2}/\mathcal{G}^{2}_{m}$ is compact since $G$ is a compact group.

\begin{proposition}
The functional $\overline{YMD}_{|\overline{\mathcal{N}}}$ satisfies the (PS) condition.
\end{proposition}
{\it Proof:} 
We will follow closely the proof of the (PS) condition for the Yang-Mills functional as in \cite{Wil}. 

Let $(A_{i},\psi_{i})$ be a (PS) sequence of $YMD_{|\mathcal{N}}$. Then 
\begin{equation}
\left\{ \begin{array}{lll}
YMD(A_{i},\psi_{i})\to c\\
\delta_{A_{i}}F_{A_{i}}-J(\psi_{i},\psi_{i})\to 0 \text{ in } \mathcal{A}^{-1}(P)\\
D_{A_{i}}\psi-K'(\psi_{i})\to 0 \text{ in } H^{-\frac{1}{2}}(\mathcal{H}).
\end{array}
\right.
\end{equation}
In particular, we have the existence of $C_{1}>0$ and $C_{2}>0$ such that
$$C_{1}\|\psi\|_{p+1}^{p+1}+\int_{M}|F_{A}|^{2}dv \leq C_{2}.$$
 Thus $\|F_{A_{i}}\|_{L^{2}}$ and $\|\psi\|_{L^{p+1}}$ are bounded. By the Uhlenbek weak compactness Theorem \cite{Uhl}, there exists a sequence of gauge transformations $(u_{i})\in \mathcal{G}^{2}(P)$ such that $u^{*}_{i}A_{i}$ is bounded in $\mathcal{A}^{1}(P)$ and weakly convergent to a connection $A_{\infty}\in A^{1}(P)$ and the convergence is strong in $\mathcal{A}^{q,0}$ for all $q\geq 1$. We will set $\tilde{A}_{i}=u^{*}A_{i}$ and $\tilde{\psi}^{i}=u^{*}_{i}\psi_{i}$, then we have that $(\tilde{A}_{i},\tilde{\psi}_{i})$ is also a (PS) sequence for $YMD_{|\mathcal{N}}$. Notice now that since $\tilde{\psi}_{i}\in \mathcal{N}$ then we have that
$$P_{\tilde{A}_{i}}^{-,0}(D_{\tilde{A}_{i}}\tilde{\psi}_{i}-K'(\tilde{\psi}_{i}))=0.$$
Hence,
$$\|\tilde{\psi}_{i}^{-}\|^{2}_{\frac{1}{2},\tilde{A}_{i}}\leq C\|\tilde{\psi}\|^{p}_{p+1}\|\tilde{\psi}^{-}\|_{\frac{1}{2},\tilde{A}_{i}}$$
Therefore $\|\tilde{\psi}_{i}^{-}\|_{\tilde{A}_{i}}$ is bounded, moreover,
\begin{align}
\|\tilde{\psi}_{i}^{+}+\tilde{\psi}_{i}^{-}\|_{\frac{1}{2},A_{i}}^{2}&=\int_{M}\langle D_{\tilde{A}_{i}}\tilde{\psi}_{i},\tilde{\psi}_{i}\rangle dv +2\|\tilde{\psi}_{i}^{-}\|^{2}_{\frac{1}{2},\tilde{A}_{i}}\notag \\
&\leq \int_{M}\langle K'(\tilde{\psi}_{i}),\tilde{\psi}_{i} \rangle +C\notag \\
&\leq C(\|\tilde{\psi}_{i}\|_{p+1}^{p+1}+1).\notag
\end{align}
Also, since $H^{\frac{1}{2},0}_{A}$ is finite dimensional, then all the norms are equivalent therefore
$$\|\tilde{\psi}_{0}\|_{\frac{1}{2},\tilde{A}_{i}}^{2}\leq C\|\tilde{\psi}_{0}\|_{p+1}.$$
Therefore, $\|\tilde{\psi}_{i}\|_{\frac{1}{2},\tilde{A}_{i}}$ is bounded and since $\tilde{A}_{i}$ is bounded in $\mathcal{A}^{1}(P)$. Using Lemma 4.2, we have that $\|\tilde{\psi}_{i}\|_{\frac{1}{2},A_{\infty}}$ is bounded. So we can extract a weakly convergent subsequence of $\tilde{\psi}_{i}$ that converges to $\psi_{\infty}$ weakly in $H^{\frac{1}{2}}(\mathcal{H})$ and strongly in $L^{q}$ for all $q<4$. Since $(\tilde{A}_{i},\tilde{\psi}_{i})$ is a (PS) sequence, we have also that
$$D_{\tilde{A}_{i}}=K'(\tilde{\psi}_{i})+o(1).$$
Using the strong convergence of $\tilde{A}_{i}$ in $L^{q}$ for all $q>1$, we deduce that
$$D_{\tilde{A}_{\infty}}\psi_{\infty}=K'(\psi_{\infty}).$$ 
Similarly,
$$\delta_{A_{\infty}}F_{A_{\infty}}=J(\psi_{\infty},\psi_{\infty}).$$
In particular, we can assume from the regularity result in Lemma 5.1 in the Appendix below, that $A_{\infty}$ and $\psi_{\infty}$ are classical solutions. So we write
$$D_{A_{\infty}}(\psi_{\infty}-\tilde{\psi}_{i})=D_{A_{\infty}}\psi_{\infty}-D_{\tilde{A}_{i}}\tilde{\psi}_{i}+R(A_{\infty}-\tilde{A}_{i})\tilde{\psi}_{i},$$
where $R(A)\psi=\sigma_{*}(A^{\alpha})e_{\alpha}\cdot \psi$ is a linear expression in $A$. Notice now that since $\tilde{\psi}_{i}$ converges strongly in $L^{2}$ and $\tilde{A}_{i}$ converges strongly in $L^{4}$, we have that $R(A_{\infty}-\tilde{A}_{i})\tilde{\psi}_{i}$ converges strongly to zero in $L^{\frac{4}{3}}(\mathcal{H})\hookrightarrow H^{-\frac{1}{2}}(\mathcal{H})$, also
\begin{align}
D_{A_{\infty}}\psi_{\infty}-D_{\tilde{A}_{i}}\tilde{\psi}_{i}&=D_{A_{\infty}}\psi_{\infty}-K'(\psi_{\infty})-(D_{\tilde{A}_{i}}\tilde{\psi}_{i}-K'(\tilde{\psi}_{i}))+K'(\psi_{\infty})-K'(\tilde{\psi}_{i})\notag\\
&=-(D_{\tilde{A}_{i}}\tilde{\psi}_{i}-K'(\tilde{\psi}_{i}))+K'(\psi_{\infty})-K'(\tilde{\psi}_{i}).
\end{align}
Since $p+1<4$, we have that $K'(\tilde{\psi}_{i})$ converges strongly in $L^{q}$ for $q<\frac{4}{p}$, then $K'(\psi_{\infty})-K'(\tilde{\psi}_{i})$ converges strongly to zero in $L^{\frac{4}{3}}\hookrightarrow H^{-\frac{1}{2}}(\mathcal{H})$. Also, since $(\tilde{A}_{i},\tilde{\psi}_{i})$ is a (PS) sequence, we have that $D_{\tilde{A}_{i}}\tilde{\psi}_{i}-K'(\tilde{\psi}_{i})$ converges strongly to zero in $H^{-\frac{1}{2}}(\mathcal{H})$. Taking $\lambda \in \R$ not  a spectral value of $D_{A_{\infty}}$, we have that
$$(D_{A_{\infty}}-\lambda)(\psi_{\infty}-\tilde{\psi}_{i})\to 0 \text{ in } H^{-\frac{1}{2}}(\mathcal{H}).$$
So by elliptic regularity of the Dirac operator, we have that $\psi_{\infty}-\tilde{\psi}_{i}\to 0$ in $H^{\frac{1}{2}}(\mathcal{H})$.

Now, using a Coulomb gauge around $A_{\infty}$ we can assume that 
$$\delta_{A_{\infty}}(\tilde{A}_{i}-A_{\infty})=0.$$
Setting $\tau_{i}=\tilde{A}_{i}-A_{\infty}$, we have that $\delta_{A_{\infty}}(\tau_{i})=0$ and
$$\Delta_{\infty} \tau_{i}=\delta_{\tilde{A}_{i}}F_{\tilde{A_{i}}}-\delta_{A_{\infty}}F_{A_{\infty}}+Q(\tau_{i}),$$
where $$Q(\tau_{i})=\frac{1}{2}\delta_{A_{\infty}}[\tau_{i},\tau_{i}]-*[\tau_{i},*(F_{A_{\infty}}+d_{A_{\infty}}\tau_{i}+\frac{1}{2}[\tau_{i},\tau_{i}])].$$
We can see that since $\tau_{i}$ converges weakly to zero in $\mathcal{A}^{1}(P)$, that $Q$ converges strongly to zero in $\mathcal{A}^{-1}(P)$. Also, we have that
\begin{align}
\delta_{\tilde{A}_{i}}F_{\tilde{A_{i}}}-\delta_{A_{\infty}}F_{A_{\infty}}&=\delta_{\tilde{A}_{i}}F_{\tilde{A_{i}}}-J(\tilde{\psi}_{i},\tilde{\psi}_{i})-(\delta_{A_{\infty}}F_{A_{\infty}}-J(\psi_{\infty},\psi_{\infty}))\notag \\
&\quad+J(\tilde{\psi}_{i},\tilde{\psi}_{i})-J(\psi_{\infty},\psi_{\infty})\notag\\
&=\delta_{\tilde{A}_{i}}F_{\tilde{A_{i}}}-J(\tilde{\psi}_{i},\tilde{\psi}_{i})+J(\tilde{\psi}_{i},\tilde{\psi}_{i})-J(\psi_{\infty},\psi_{\infty}).\notag
\end{align}
 Again, since $(\tilde{A}_{i},\tilde{\psi}_{i})$ is a (PS) sequence, we have that $ \delta_{\tilde{A}_{i}}F_{\tilde{A_{i}}}-J(\tilde{\psi}_{i},\tilde{\psi}_{i})$ converges strongly to zero in $\mathcal{A}^{-1}(P)$ and since $\tilde{\psi}_{i}$ converges strongly in $L^{q}$ for $q<4$, we have that $J(\tilde{\psi}_{i},\tilde{\psi}_{i})-J(\psi_{\infty},\psi_{\infty})$ converges strongly to zero in $L^{2}\hookrightarrow \mathcal{A}^{-1}(P)$, hence $\Delta_{A_{\infty}}\tau_{i}$ converges strongly to zero in $\mathcal{A}^{-1}(P)$. Again using the compactness of the operator $\Delta_{A_{\infty}}+1$, we have the strong convergence of $\tau_{i}$ to zero in $\mathcal{A}^{1}(P)$, which finishes the proof since $\mathcal{G}^{2}(P)/\mathcal{G}_{m}^{2}$ is compact.
\hfill$\Box$

\begin{proposition}
The functional $YMD$ has infinitely many non-gauge equivalent, non-trivial solutions.
\end{proposition}
{\it Proof:}

To prove the following, it is enough to show that $\overline{\mathcal{N}}$ has infinite genus. For that, we fix $A\in \mathcal{A}^{1}(P)$ and we consider the map $Z: S^{+}_{A}\to \mathcal{N}$ defined by
$$Z(\psi)=(A,T_{A}(\psi))$$
The set $S_{A}^{+}$ is invariant under the action of $\Z_{2}$, that is $-S^{+}_{A}=S^{+}_{A}$. Moreover, we have that 
$$T_{A}(-\psi)=-T_{A}(\psi).$$
Thus the map $Z$, is equivariant, and since $\gamma(S^{+}_{A})=\infty$, we have that $\gamma(\overline{\mathcal{N}})=\infty$, therefore, by Proposition 4.6,  if we denote by $\mathcal{B}_{k}$ the collection of sets $B\subset \overline{\mathcal{N}}$ such that $\gamma(B)\geq k$, have that the values
$$c_{k}=\inf_{B\in \mathcal{B}_{k}}\max_{B}\overline{YMD}$$

are critial values of $\overline{YMD}$, which finishes the proof.
\hfill $\Box$

\section{Appendix: Regularity}

We consider now a weak solution of the system
\begin{equation}\label{eqy}
\left\{\begin{array}{ll}
\delta_{A}F_{A}=J(\psi,\psi)\\
D_{A}\psi=K'(\psi).
\end{array}
\right.
\end{equation}
\begin{lemma}
If $(A,\psi)$ is a solution to $(\ref{eqy})$ then there exists $u\in \mathcal{G}^{2,2}$ such that $(u^{*}A,u^{*}\psi)\in C^{2,\alpha}\times C^{1,\alpha}$.
\end{lemma}
{\it Proof:}

We will place our selves in a Coulomb gauge with respect to a smooth connection $A_{0}$ close to $A$ in the $\mathcal{A}^{1}(P)$ norm. We will assume without loss of generality that $A_{0}=0$. That is we will replace $(A,\psi)$ by $(u^{*}A,,u^{*}\psi)$ so that  we have
$$\delta_{0}u^{*}A=0$$
We will disregard from now on the action of $u$. Thus, we have that
$$\delta_{0}(F_{A})-*[A,*F_{A}]=J(\psi,\psi).$$
Notice then that we have
$$\Delta_{0}A=-\frac{1}{2}\delta_{0}[A,A]+*[A,*F_{A}]+J(\psi,\psi)$$
and 
$$D_{0}\psi=K'(\psi)-R(A)\psi.$$
Since $A\in \mathcal{A}^{1}(P)$, we have that $A\in L^{q}$ for all $q>1$ and $dA\in L^{2}$, hence $*[A,*F_{A}]\in L^{q}$ for all $q<2$ and similarly for $\delta_{0}[A,A]$. Also since $\psi \in H^{\frac{1}{2}}(\mathcal{H})$, then we have that $\psi \in L^{r}$ for all $1\leq r \leq 4$. Therefore we have that $\Delta_{0}A\in L^{q}$ for all $q<2$, by classical elliptic regularity, we have that $A\in \mathcal{A}^{2,q}(P)$ hence $A\in C^{0,\alpha}$.

On the other hand, we have that $K'(\psi)\in L^{\frac{4}{p}}$ thus $D_{0}\psi \in L^{\frac{4}{p}}$ again by elliptic regularity, we have that $\psi\in W^{1,\frac{4}{p}}(\mathcal{H})$. Here, we have different cases.

\textbf{Case 1:} If $1<p\leq 2$.\\
Then $\psi \in L^{r}$ for all $r>1$. Hence, $D_{0}\psi \in L^{r}$ for all $r>1$ so $\psi\in C^{0,\alpha}$, iterating again using Schauder's estimates, we have that $\psi \in C^{1,\alpha}$. 

\textbf{Case 2:} If $3>p>2$. \\
Then $\psi \in L^{r}$ and $r=\frac{4}{p-2}>4$, so  by classical boot-strap argument, we have that $\psi \in C^{0,\alpha}$, once again using Schauder's estimates we have that $\psi \in C^{1,\alpha}$.
 
Now, we go back to $A$. Notice that $\Delta_{0}A\in \mathcal{A}^{1,q}(P)$ for $q<2$. Thus $A\in \mathcal{A}^{3,q}(P)$, hence, $A\in C^{2,\alpha}$.
\hfill$\Box$.

\end{document}